\documentclass[11pt, a4paper]{article}
\usepackage{authblk}
\usepackage{amsthm}  
\usepackage{amssymb}
\usepackage{amsmath}
\usepackage{graphicx}
\usepackage{hyperref}  

\newtheorem{theorem}{Theorem}[section]
\newtheorem{theoremN}[theorem]{Theorem}
\numberwithin{theorem}{section}
\newtheorem{definitionN}[theorem]{Definition} 
\newtheorem{lemmaN}[theorem]{Lemma}
\newtheorem{corollaryN}[theorem]{Corollary}
\newtheorem{notationN}[theorem]{Notation}
\newtheorem{terminologyN}[theorem]{Terminology}
\newtheorem{remarkN}[theorem]{Remark}
\newtheorem{observationN}[theorem]{Observation}
\newtheorem{proofOfResultN}[theorem]{}
\newenvironment{proofN}{\trivlist \item[\hskip \labelsep{\bf Proof.}]}{
$\blacksquare$ \endtrivlist   }

\begin{document}

\title{\bf Desargues theorem, its configurations, and the solution to a long-standing enumeration problem
}


\author[1]{Aiden A. Bruen
  \thanks{The first author gratefully acknowledges the financial support of the National Sciences and Engineering Research Council of Canada}
}
\author[2]{Trevor C. Bruen}
\author[3]{James M. McQuillan}

\affil[1]{School of Mathematics and Statistics,
  Carleton University, 
  1125~Colonel By Drive,
  Ottawa, Ontario, K1S 5B6, Canada\\
  {\tt abruen@math.carleton.ca}
}
\affil[2]{Universit\'e de Sherbrooke,
  Sherbrooke, Quebec, J1K 2R1, Canada\\
  {\tt tbruen@gmail.com}
}
\affil[3]{School of Computer Sciences, 
  Western Illinois University,
  1~University Circle,
  Macomb, IL, 61455, USA\\
 {\tt jm-mcquillan@wiu.edu}
}

\date{}
\maketitle

\begin{abstract}
  We solve a long-standing problem by enumerating the number
  of non-degenerate Desargues configurations.
  We extend the result to the more difficult case
  involving Desargues blockline structures
  in Section~\ref{section:DesarguesBlockline}.

  A transparent proof of Desargues theorem in the plane 
  and in space
  is presented as a by-product of our methods.
\end{abstract}


%
%
%
\thispagestyle{empty}

\section{\ Introduction, Background.\label{introduction}}
The celebrated theorems of Pappus and Desargues
are the fundamental building blocks in 
the axiomatic development of incidence and projective
geometry. Hilbert's work showed that
a projective plane $\pi$ over a division 
ring $F$ is equivalent to one in which
Desargues theorem holds,
which is in turn equivalent to the assumption
that $\pi$ is embedded in a 3-dimensional projective space.

The Hessenberg theorem shows that the ancient theorem of
Pappus implies the Desargues theorem.  From this
it follows that a projective plane over a commutative
division ring $F$, i.e. a field $F$, is
equivalent to one in which the Pappus theorem holds.
Thus, since all finite division rings are fields,
we have that Pappus and Desargues are equivalent
for finite projective planes.
For a discussion of non-Desarguesian planes see Lorimer~\cite{lorimer}.
The Desargues configuration still plays a fundamental role
when studying the collineation of such planes.

In~\cite{rota} page 145, the author, referring also to 
Baker~\cite{baker} writes as follows:
``After an argument that runs well over one hundred pages,
Baker shows that beneath the statement of Desargues' theorem,
another far more interesting geometric structure lies {\em concealed}.
This structure is nowadays called the Desargues configuration.''

As an aside,
the Desargues configuration makes its appearance also in many
areas of combinatorics and graph theory.
In the paper on Colouring Problems by W.T.~Tutte
(see~\cite{tutte}),
the author discusses connections between 
graph theory and the geometrical approach to the
Four Colour Problem first explored by
O. Veblen~\cite{veblen}
in 1912.

In a 3-dimensional projective space over the binary
field a {\em tangential 2-block} is a set $S$ of points
with the following two properties:
\begin{itemize}
  \item[(a)]
    Each line in space contains at least one point of $S$;
  \item[(b)]
    At each point $P$ of $S$ there exists a tangent line
    to $S$ at $P$, i.e., a line meeting $S$ only in $P$.
\end{itemize}
A remarkable fact is that, apart from a plane in the
space, i.e., a Fano plane, the only other tangential
2-block is a spatial (i.e., a non-planar) Desargues configuration.

Generalizing the definition of tangential 2-block
to projective spaces of arbitrary dimension,
Tutte announces in \cite{tutte} his famous tangential 2-block conjecture
to the effect that the only tangential 2-blocks existing
in projective space over the binary field
correspond to the Fano plane,
the Desargues configuration and the 5-dimensional
Petersen block.
The conjecture, as a special case,
implies the Four Colour Theorem for graphs.
The conjecture has been verified for
several dimensions
but remains unsolved as of this writing.

Enumerating the number of Desargues configurations
in the plane is an old problem dating back to 1970-71
and the work in~\cite{thas}.  There the authors
restrict attention to the special case where the characteristic
of the underlying field is~2 or~3 and no line of the plane
contains as many as 4 points of the 10 points of the configuration.
This is tantamount to saying that no point of the
configuration is self-conjugate.
Later on we show in Theorem~\ref{theorem:blocklineAtMost4Pts}
that no line
contains more than 4 points of the configuration.
For further discussion on conjugacy we refer to 
our paper~\cite{bruenMcQGeomConfigs}.

Even in the special case discussed in~\cite{thas}
the proof is lengthy and complicated.
In this paper, we produce a complete solution.
We show in Section~\ref{section:easierApproach}
how the standard approach doesn't work
given the many possibilities that must be considered
and the impact of the characteristic of the field.
A crucial factor here is that the type of configuration
depends on the number of self-conjugate points,
which in turn depends on the characteristic of the field.
We also enumerate the 3-dimensional Desargues configurations,
providing two proofs of the result.

%
%
%
%
%
\section{\ Projective planes and spaces, Desargues configurations. \label{section:DesarguesConfigurations}}
For basic definitions, we refer to 
Coxeter~\cite{coxeter,coxeter2},
Pedoe~\cite{pedoe},
Hirschfeld~\cite{hirschfeld},
Hartshorne~\cite{hartshorne},
Veblen and Young~\cite{veblenYoungVol1},
or Todd~\cite{todd}.
For some interesting related work, we mention
articles by Conway and Ryba~\cite{conway,conway2},
Crannell and Douglas~\cite{crannell}
and Lord's book~\cite{lord}.
A recent paper of the authors~\cite{bruenMcQGeomConfigs}
relates to the material in Section~\ref{section:DesarguesBlockline}.

Much of the paper can be visualized 
in projective space over the reals.  Here is a quick overview
from the synthetic or algebraic point of view.
The Euclidean plane can be embedded
in, or extended to the real
projective plane $PG(2,\mathbb{R})$
as follows.  We adjoin an ``infinite point''
(or ``slope point'') to each of
a given parallel class of lines.
Different parallel classes have different slope points.
We decree that the newly created
infinite points lie on an
``infinite line'' or ``line at infinity''.
The projective plane now has the property that
two distinct points lie on a unique line
and two distinct lines meet in a unique point.

There is also the more algebraic approach
as follows.  Start with a 3-dimensional vector
space $V$ over the reals.
Define points to be all the
1-dimensional subspaces of $V$
and the lines to be all the 2-dimensional
subspaces of $V$.
The resulting structure is isomorphic to
$PG(2,\mathbb{R})$ as described above
as shown in Hartshorne~\cite{hartshorne}.

Similarly, Euclidean 3-space can be embedded
in the projective space  \hfil\break
$PG(3,\allowbreak\mathbb{R})$ by adjoining
a ``plane at infinity''.
The algebraic model then comes from a 4-dimensional
real vector space.

Initially we work in $\Sigma=PG(n,F)$
the $n$-dimensional projective space based on the $(n+1)-$dimensional
vector space over the skew field $F$.
The celebrated theorem of Desargues is as follows.
\begin{figure}[!ht]
    \center
    (a)
    \includegraphics[height=7.9cm]{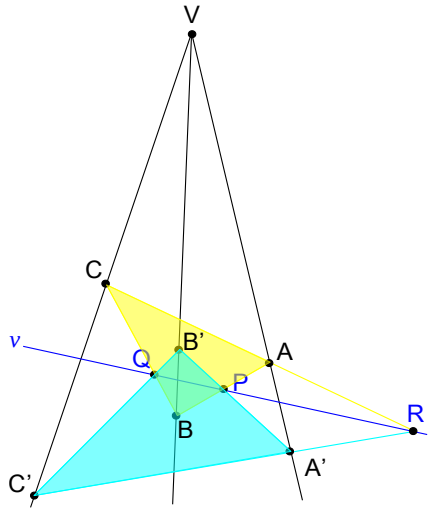}\\
    (b)
    \includegraphics[height=7.9cm]{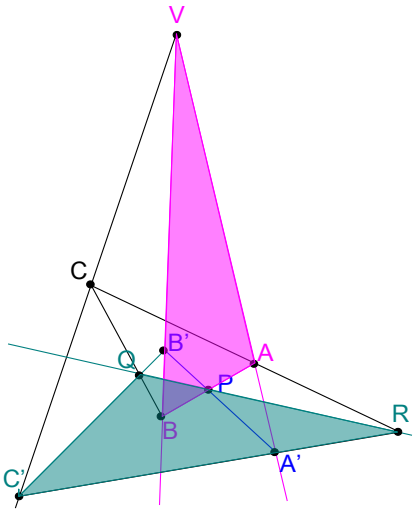}
    \caption{A Desargues configuration.
      (a) if we choose the point $V$ for the vertex,
      then $ABC$ and $A'B'C'$ are in perspective from $V$
      with Desargues axis $v$ containing $P,Q,R$.
      (b) if we choose the point $C$ for the vertex,
      then $QRC'$ and $BAV$ are in perspective from $C$
      with Desargues axis containing $A',B',P$.
      The block $A'B'P$ is the polar of $C$.
    }
    \label{fig:desConfig}
\end{figure}
\begin{theoremN}
  Let $ABC$ and $A'B'C'$ be two triangles in $\Sigma$
  such that the lines $AA'$, $BB'$ and $CC'$ pass through
  a point $V$ of $\Sigma$.
  Then the 3 intersection points $P,Q,R$ of the lines
  $AB$ with $A'B'$, $BC$ with $B'C'$ and
  $CA$ with $C'A'$ lie on a line $v$.
  \label{thm:DesarguesTheorem}
\end{theoremN}
We refer to Figure~\ref{fig:desConfig}(a).
Note that the two triangles either lie in a plane
or are contained in a 3-dimensional space 
so we may assume that $n=3$.
The 2 triangles are said to be
{\em in perspective} with center $V$.
For reasons explained below and in Section~\ref{section:DesPlanePolarities},
the line $v$ is called
the {\em Desargues axis} or {\em polar} of $V$.

The two triangles yield a {\em Desargues Configuration} $D$
with 10 points and 10 triples of collinear points called
{\em blocks}, as follows.
\begin{itemize}
  \item[]
    Points of $D$: $V,A,B,C,A',B',C',P,Q,R$
  \item[]
    Blocks of $D$: $\{VAA'\}$, $\{VBB'\}$, $\{VCC'\}$,
    $\{APB\}$, $\{A'PB'\}$, \\
    \hphantom{Blocks of $D$: }
    $\{ARC\}$, $\{A'RC'\}$,
    $\{BQC\}$, $\{B'QC'\}$, $\{PQR\}$.
\end{itemize}
A line of $\Sigma$ containing a block of $D$ is called a
{\em blockline}.
Each point of $D$ is contained in exactly 3 blocks of $D$
and each block of $D$ contains exactly 3 points of $D$.

{\em Throughout this paper, we assume 
that $D$ is non-degenerate, i.e., that $D$
has 10 distinct points and 10 distinct blocklines.}

The configuration $D$ is an example of a $(3,3)$-configuration
in the terminology of Pedoe~\cite[p. 25]{pedoe},
or a $10_3$-configuration in the terminology of
Coxeter~\cite[p. 26]{coxeter}.

As with $V$, any of the 10 points of $D$
can be regarded as the vertex of perspective of
exactly two triangles formed from 
the points and blocklines of $D$.
For example, if we choose the point $C$
for vertex then, using the 3 blocklines on $C$,
we get the 2 triangles in perspective from $C$,
namely, the triangles $QRC'$ and $BAV$.
The Desargues axis for the vertex $C$ is the blockline
containing the block $\{A'B'P\}$.
See Figure~\ref{fig:desConfig}(b).

%
%
%
\section{\ Planar Desargues configurations and 5-compressors.\label{sect:planarDesAnd5Pts}}
Here we work with {\em planar Desargues configurations}
$D$ where the 2 triangles $ABC$ and $A'B'C'$
lie in the same plane lying in $\Sigma=PG(3,F)$,
with $F$ any field.

\begin{definitionN}
  A {\em 5-compressor} or a {\em 5-point} or a 
{\em simplex $\mathcal{S}$} in $\Sigma$ is a set $\mathcal{S}$
of 5 points with no 4 coplanar.
\end{definitionN}
In what follows, the symbol
${<}P_1,P_2,P_3{>}$ represents
the plane containing the triangle 
$P_1P_2P_3$.

We now develop the important connection between 5-compressors
and Desargues configurations.

\begin{theoremN}
  Let $\mathcal{S}=\{P_1,P_2,P_3,P_4,P_5\}$
  be a 5-compressor,
  i.e., a set of 5 points in $\Sigma=PG(3,F)$
  with no~4 of its points coplanar.
  Let $\pi$ be any plane of $\Sigma$ containing no point of $\mathcal{S}$.
  The points $(ij)$, where $(ij)$ is the point of intersection
  of the line $P_iP_j$ with $\pi$, yield~10
  distinct points in $\pi$.
  The plane containing $P_i,P_j,P_k$ intersects $\pi$
  in a line containing the block $[ijk]$ which denotes
  the set of~3 collinear points $(ij)$, $(jk)$ and $(ik)$.
  The 10 points $(ij)$ and the 10 blocks such as $[ijk]$
  form a Desargues $(3,3)$-configuration in~$\pi$.
  Conversely, any Desargues configuration
  arises from
  a 5-compressor in $\Sigma$
  using the above construction.
  \label{5Point}
\end{theoremN}
\begin{proofN}
  The assumption on $\mathcal{S}$ guarantees that
  no~3 points $P_i,P_j,P_k$ of $\mathcal{S}$ are collinear.
  If the line $P_iP_j$ were to also meet $\pi$ in the point
  $P_uP_v\cap\pi$,
  the sets
  $\{u,v\}$ and $\{i,j\}$ being disjoint
  would imply that $\{P_i,P_j,P_u,P_v\}$
  is a set of~4 points lying in a plane.
  This contradicts the definition of~$\mathcal{S}$.
  Thus the 10 lines $P_iP_j$ meet $\pi$
  in ${5 \choose 2}=10$ distinct points.

  A plane containing $\{P_i,P_j,P_k\}$, denoted by ${<}P_i,P_j,P_k{>}$,
  meets $\pi$ in~3 distinct collinear points.
  Suppose that~2 planes $\pi_1=\{P_i,P_j,P_k\}$ and 
  $\pi_2=\{P_u,P_v,P_w\}$
  meet $\pi$ in a line $l$.
  Since no line contains more than 4 points of the
  configuration, it follows that
  $P_u$, say, must equal one of $P_i,P_j,P_k$, say $P_i$.
  Then $\pi_1$ and $\pi_2$ are both equal to
  the plane $\pi_3$ containing $P_u=P_i$ and $l$.
  Therefore $\pi_3$ contains either 4 or 5 points of $\mathcal{S}$,
  which is impossible.

  In this way we obtain
  ${5\choose 2}=10$ points $(ij)$
  and ${5\choose 3}=10$ blocks $[ijk]$
  corresponding to the lines  
  $P_iP_j$
  and the planes ${<}P_i,P_j,P_k{>}$.

  Points $(ab)$ and $(cd)$ lie in  one of
  the~10 blocks if and only if
  the pairs $(ab)$ and $(cd)$ share a symbol.
  Two of the 10 blocks $[abc]$ and $[uvw]$
  intersect in one of the 10 points~$(ij)$
  if and only if
  the~2 triples share a common pair.
  In this way we see that the 5-compressor $\mathcal{S}$
  yields a $(3,3)$-Desargues configuration $D$
  in $\pi$
  with 10 points and 10 lines.

  We will prove a strong form of the converse
  in Theorem~\ref{thm:two5Points}.
\end{proofN}

In light of Theorem~\ref{5Point}, we can think of
5-compessors as giving us data compression ---
high above the Desargues configuration, ``in the clouds'' ---
in that they store the information of Desargues configurations.

\begin{figure}[!ht]
    \center
    \includegraphics[height=9cm]{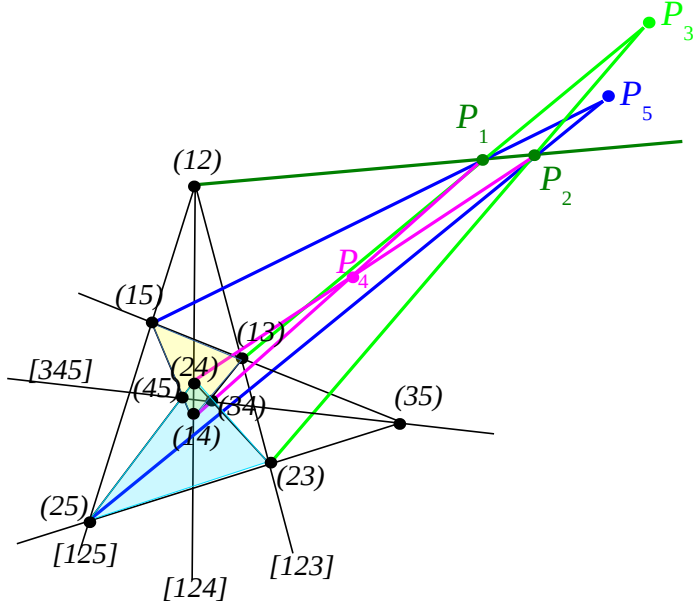}
    \caption{A 5-compressor $\mathcal{S}=\{P_1,P_2,P_3,P_4,P_5\}$
      and the Desargues configuration $D$ in $\pi$ that $\mathcal{S}$ sections.
      The Desargues configuration $D$ is shown in more detail in
      Figure~\ref{fig:alternativeLabels}.}
    \label{fig:aSimplex}
\end{figure}

\begin{figure}[!ht]
    \center
    \includegraphics[height=14cm]{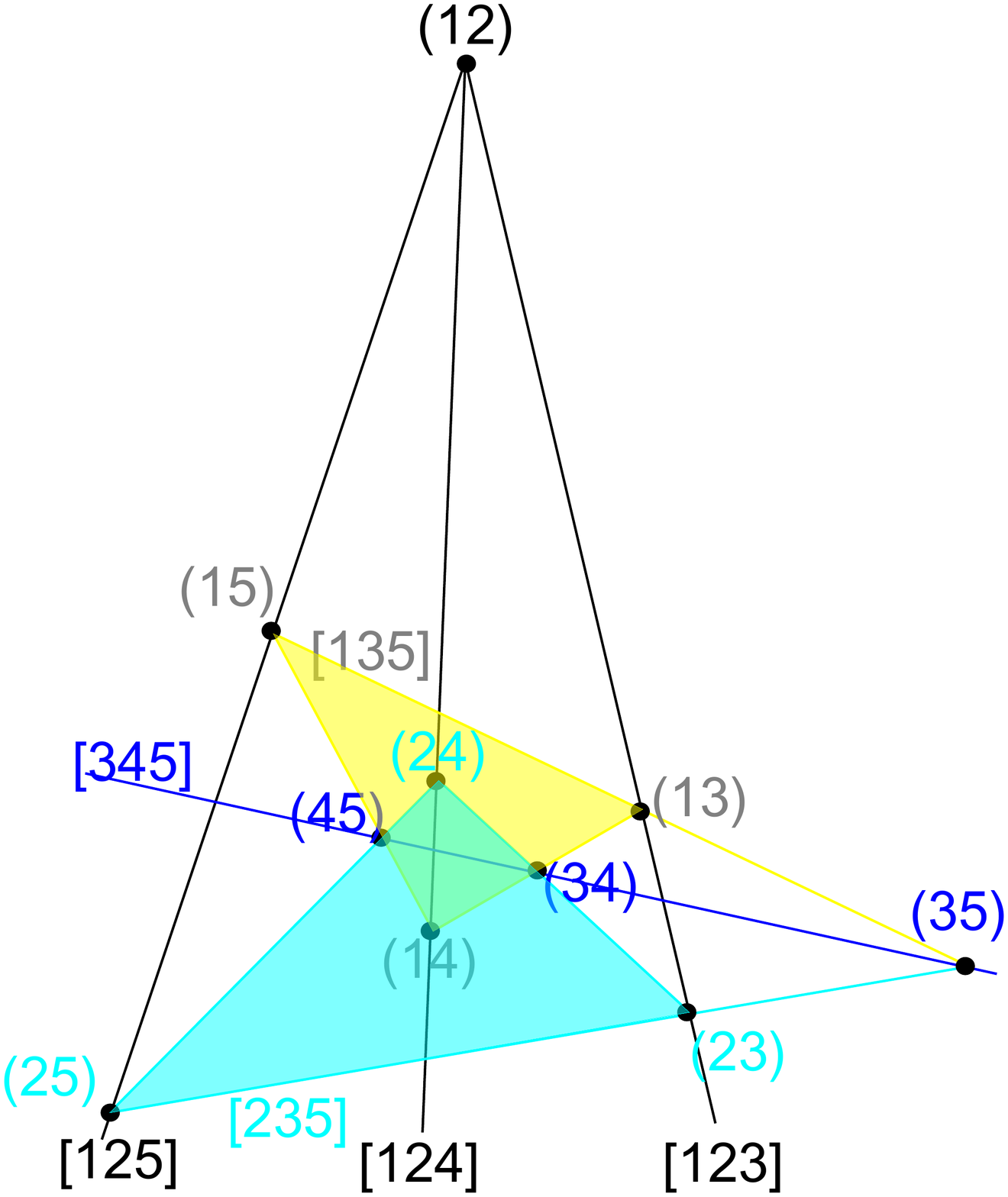}
    \caption{Alternative labels for the Desargues configuration
      in Figure~\ref{fig:desConfig}(a).  Here,
      $V=(12)$, $A=(13)$, $A'=(23)$,
      $B=(14)$, $B'=(24)$,
      $C=(15)$, $C'=(25)$,
      $P=(34)$, 
      $Q=(45)$, 
      and $R=(35)$.
      For this labeling, 
      the points of block $[123]$ are on blockline $AA'$, 
      $[124]$ is on $BB'$,
      $[125]$ is on $CC'$,
      $[134]$ is on $AB$, 
      $[234]$ is on $A'B'$,
      $[135]$ is on $AC$, 
      $[235]$ is on $A'C'$, 
      $[145]$ is on $BC$, 
      $[245]$ is on $B'C'$, and
      $[345]$ is on $PQ$.
    }
    \label{fig:alternativeLabels}
\end{figure}

We refer to Figure~\ref{fig:aSimplex} for an illustration
of a 5-compressor, as described in Theorem~\ref{5Point}.
All of the lines $P_1P_i$ and $P_2P_j$ are shown,
$i\ne 1$, $j\ne 2$.
The Desargues configuration $D$ is shown in more detail in
Figure~\ref{fig:alternativeLabels}.


\begin{terminologyN}
  We refer to the configuration~$D$
  with its 10 points and 10 blocks as above
  as the {\em section} of the 5-compressor (5-point) $\mathcal{S}$ by $\pi$,
  and we say that {\em $\mathcal{S}$ sections $D$}.
  \label{section5PointPi}
\end{terminologyN}
This is consistent with the terminology in~\cite[p. 34]{veblenYoungVol1}.

\begin{figure}[!ht]
    \center
    (a) 
    \includegraphics[height=8.7cm]{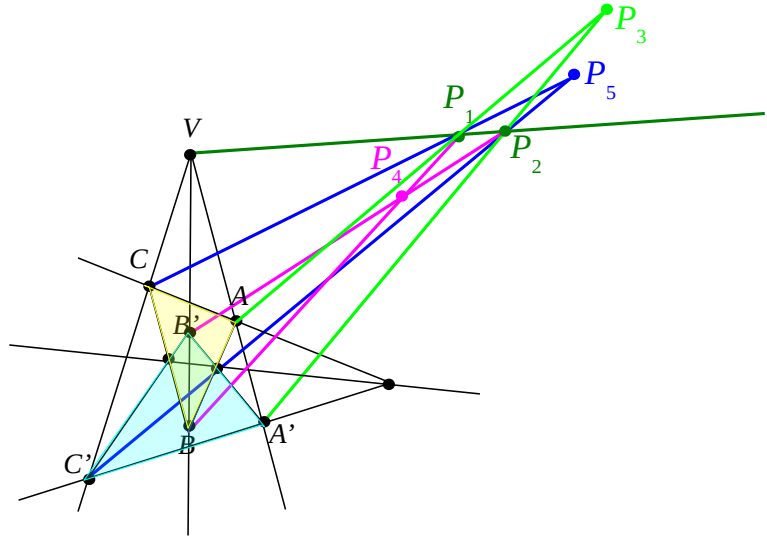}
    \\
    (b) 
    \includegraphics[height=8.7cm]{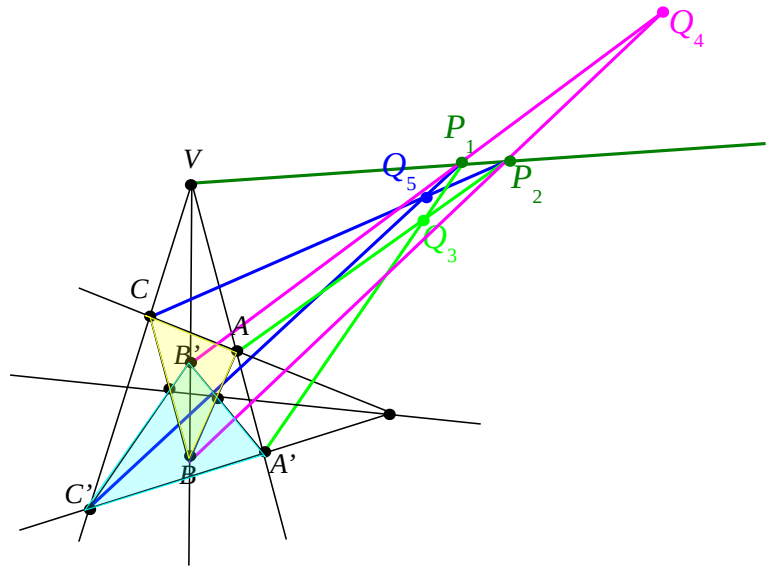}
    \caption{The 5-compressors $\mathcal{S}_1$ in (a),
      and $\mathcal{S}_2$ in (b),
      section to $D$ for 
      Theorem~\ref{thm:ExistsTwo5Pts}.
      Alternatives labels for the
      points and blocks
      of $D$ corresponding to $\mathcal{S}_1$
      are shown in Figure~\ref{fig:alternativeLabels}.
    }
    \label{fig:twoSimplexes}
\end{figure}

\begin{theoremN}
  Let $D$ be any (non-degenerate) Desargues configuration
  in the plane $\pi=PG(2,F)$, with $F$ any field.
  Let $V$ be any point of $D$.
  Let $P_1,P_2$ be two fixed points in space not in $\pi$
  such that $P_1$, $P_2$ and $V$ are collinear.
  Then there exists two 5-compressors,
  $\mathcal{S}_1=\{P_1,P_2,P_3,P_4,P_5\}$
  and
  $\mathcal{S}_2=\{P_1,P_2,Q_3,Q_4,Q_5\}$
  which section to $D$,
  where $\mathcal{S}_1$, $\mathcal{S}_2$ contain
  no point of $\pi$.
  \label{thm:ExistsTwo5Pts}
\end{theoremN}
\begin{proofN}
  Let $D$ be labelled as indicated in
  Figure~\ref{fig:desConfig}(a).
  Let $P_3$ denote $P_1A\cap P_2A'$,
  i.e., $P_3$ is the point of intersection of the lines
  $P_1A$ and $P_2A'$ joining $P_1$ to $A$ and $P_2$ to $A'$
  respectively.
  (Alternatively let 
  $Q_3=P_1A'\cap P_2A$.)
  Note that $Q_3$ cannot equal $P_3$ because for example
  $P_1$ is not in the plane $\pi$.)
  Let $P_4=P_1B\cap P_2B'$ and $P_5=P_1C\cap P_2C'$.
  (Alternatively,
  let $Q_4=P_1B'\cap P_2B$ and $Q_5=P_1C'\cap P_2C$.
  Similar to the above, $Q_4$ cannot equal $P_4$ and $Q_5$ cannot equal $P_5$.)
  As in the proof of Theorem~\ref{5Point}, 
  using either $\mathcal{S}_1$ or $\mathcal{S}_2$ we obtain the 10 points
  $(ij)$ and the 10 blocks $[ijk]$ of a Desargues configuration.

  Here, using $\mathcal{S}_1$ we get
  $V=(12)$, $A=(13)$, $A'=(23)$, $B=(14)$,
  $B'=(24)$, $C=(15)$, and $C'=(25)$.
  We define $Q$ as the intersection of the blocks
  $CB$ and $C'B'$. 
  The block $CB$ contains $C=(15)$ and $B=(14)$
  so the block $CB$ contains $(15)$, $(14)$ and $(45)$.
  The block $C'B'$ contains $(25)$, $(24)$ and $(45)$. 
  Thus $Q=(45)$.
  Similarly $P$ is $(34)$ and $R=(35)$.
  Then $P,R,Q$ are collinear on the blockline containing
  $(34),(35),(45)$, i.e., the block $[345]$.
  We refer to
  Figures~\ref{fig:alternativeLabels}
  and~\ref{fig:twoSimplexes}.

  In a similar way we can calculate that all of the 
  10 blocks in $D$ are of the form $[ijk]$.
  (Using $\mathcal{S}_2$ we get a different labelling
  obtained by interchanging $ABC$ with $A'B'C'$.)
  This completes the proof of Theorem~\ref{thm:ExistsTwo5Pts}.
\end{proofN}

%
%
%
%
%
\section{\ Desargues theorem in the plane and polarities.\label{section:DesPlanePolarities}}
Using the notation of Section~\ref{section:DesarguesConfigurations},
we now use Theorem~\ref{thm:ExistsTwo5Pts}
to give a new transparent proof of the Desargues theorem in the plane.
Unlike several standard proofs we do not explicitly rely on the
Desargues theorem in space.
Instead we use the fact that two planes intersect in a line.

\bigskip\bigskip
\begin{proofOfResultN}
  Proof of Theorem~\ref{thm:DesarguesTheorem},
  the theorem of Desargues in the plane.
\end{proofOfResultN}
\begin{proofN}
  Construct the 5-compressor $\mathcal{S}_1$ sectioning the
  given Desargues configuration in the plane $\pi$
  as in
  Theorem~\ref{thm:ExistsTwo5Pts}.
  The intersections of the pairs of corresponding sides are
  $P,Q,R$ which are collinear since
  the 3 points lie on the intersection
  of two planes, namely the plane $\pi$ and the plane
  ${<}P_3,P_4,P_5{>}$.
\end{proofN}

\begin{remarkN}
  The diagram in 
  Figure~\ref{fig:twoSimplexes}
  and
  Figure~\ref{fig:aSimplex} 
  is precisely the same as that in the standard
  proof of the theorem.  See Coxeter~\cite{coxeter}. 
\end{remarkN}

\begin{figure}[!ht]
    \center
    \includegraphics[height=10cm]{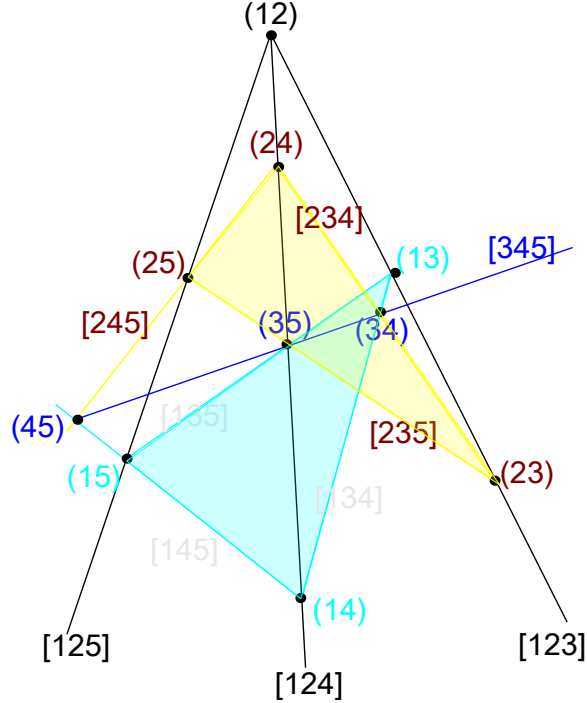}
    \caption{The triangles $(13)(14)(15)$ and $(23)(24)(25)$
      are in perspective from $(12)$.  The point $(35)$ is 
      self-conjugate; it is on the blockline containing $[124]$.
    \label{fig:twoSCPtsBlockline}}
\end{figure}

Theorem~\ref{thm:ExistsTwo5Pts}
provides us with a convenient
notation for the 10 points $(ij)$ and 
10 blocks $[ijk]$ of
a Desargues configuration $D$
as in Figure~\ref{fig:alternativeLabels}.
The mapping 
$$\phi:(ij)\longrightarrow \text{ blockline containing }[uvw]$$
and 
$$\phi: \text{blockline containing }[uvw]\longrightarrow (ij),$$
where $i,j,u,v,w$ are distinct,
has period 2,
maps points to blocklines, blocklines to points
and preserves incidences in $D$.
Thus $\phi$ is a {\em polarity} of $D$,
mapping a point to its polar blockline,
and a blockline to its pole.
If a point $P$ lies on its polar blockline $l_P$,
we say that this point $P$ is
{\em self-polar}, or that $P$ is a {\em self-conjugate (SC) point}.
Following the notation of~\cite{bruenMcQGeomConfigs},
we think of $P$ as an ``accidental'' extra point
on the blockline containing its polar block.
Dually, if a line $l$ contains its pole,
we say that $l$ is {\em self-polar} or {\em self-conjugate}.
For example, in Figure~\ref{fig:twoSCPtsBlockline},
the point $(35)$ is an SC point on the blockline 
containing $(12),(14),(24)$.

As pointed out in Section~\ref{section:DesarguesConfigurations},
each of the 10 points $(ij)$ of $D$
serves as the vertex of perspective for
two triangles in $D$.
The polar of $(12)$ is the corresponding Desargues axis
which is the blockline containing the block $[345]$.

Next, let $l$ be a blockline containing the block
$[123]$  consisting of the points $(12),(13),(23)$.
Let $l$ contain an extra point, $(ij)$.
Suppose that $i$ (or $j$)
is equal to 1 or 2 or 3.
For example, let $i=1$, $j=5$.
Then $l$ contains $(12)$, $(13)$, $(23)$, $(15)$.
Since $l$ contains $(12)$ and $(15)$, it contains
$(25)$.
Then the 5 points $V,A,A',C,C'$ lie on $l$
so that two of the blocklines are equal.
This contradicts the fact that $D$ is non-degenerate.

We summarize as follows.
\begin{theoremN}
  A blockline containing the block $[ijk]$
  can contain at most one additional point of $D$,
  namely the point $(uv)$, where $u,v,i,j,k$
  are 5 different symbols.
  In this case, the point $(uv)$ lies on the 
  Desargues axis $[ijk]$ of the two triangles
  in $D$ that are in perspective from $(uv)$.
  Thus $(uv)$ lies on its polar line,
  i.e., it is self-polar ($=$self-conjugate).
  For a detailed discussion of the polarity,
  see~\cite{bruenMcQGeomConfigs}
  and Section~\ref{section:DesarguesBlockline}.
  \label{theorem:AtMostOneAddlPt}
\end{theoremN}

It is perhaps worth restating
Theorem~\ref{theorem:AtMostOneAddlPt}
since the result is so useful.
\begin{theoremN}
  Any blockline contains at most four points of
  the Desargues configuration.
  \label{theorem:blocklineAtMost4Pts}
\end{theoremN}

%
%
%
%
\section{\ 5-Compressors sectioning to $D$.\label{5PointsSectioningToD}}
Let $V$ be a fixed point in $PG(2,F)$ and let
$D$ be a (non-degenerate) Desargues configuration containing $V$.
Let $P_1,P_2$ be two fixed points in space
not lying in $\pi=PG(2,F)$ such that $P_1$, $P_2$ and $V$
are collinear.  

The following result is a strengthening of
Theorem~\ref{thm:ExistsTwo5Pts}.
We use the notation in Figure~\ref{fig:desConfig} (a).

\begin{theoremN}
  There exist exactly two 5-compressors $\mathcal{S}_1,\mathcal{S}_2$
  with no point of either lying in $\pi$ which
  contain $P_1$ and $P_2$ and whose section is $D$.
  \label{thm:two5Points}
\end{theoremN}
\begin{proofN}
  From 
  Theorem~\ref{thm:ExistsTwo5Pts}
  there are at least two 5-compressors
  sectioning to $D$.
  Let the 5-compressor $\mathcal{S}=\{P_1,P_2,P_3,P_4,P_5\}$
  section to $D$.
  There are 3 blocks, namely $VAA'$, $VBB'$ and $VCC'$
  containing $V$.  We can assume that the 3 planes
  ${<}P_1,P_2,P_3{>}$, ${<}P_1,P_2,P_4{>}$,
  ${<}P_1,P_2,P_5{>}$ intersect $\pi$
  in the 3 blocklines containing
  $VAA'$, $VBB'$, $VCC'$ respectively.
  The blockline containing $VAA'$ might well contain
  a fourth point $X$ of $D$ with $X$ being the pole
  of the block $VAA'$.  However, neither of the
  triples $\{VAX\}$, $\{VA'X\}$ are blocks of $D$.
  It follows that we only have 2 possibilities, namely
  \begin{itemize}
    \item[(a)]
      $P_3$ is the point $P_1A\cap P_2A'$
    \item[(b)]
      $P_3$ is the point $P_1A'\cap P_2A$.
  \end{itemize}
  Assume that case (a) holds.  Examining the blockline
  containing $\{V,B,B'\}$ we again have 2 possibilities
  \begin{itemize}
    \item[(a)]
      $P_4=P_1B\cap P_2B'$
    \item[(b)]
      $P_4=P_1B'\cap P_2B$.
  \end{itemize}
  Assume case (b) here.  We already have 
  $P_3=P_1A\cap P_2A'$ so that $P_1P_3A$ are collinear.
  We now assume that $P_4=P_1B'\cap P_2B$
  so that $P_1P_4B'$ are collinear.
  By hypothesis the points of intersection of the lines
  $P_1P_3$, $P_1P_4$ and $P_3P_4$ with $\pi$
  are the 3 points of a block which lies on the blockline
  formed by the intersection of ${<}P_1,P_3,P_4{>}$
  with $\pi$.
  Now $P_1P_3$ meets $\pi$ in $A$ and $P_1P_4$ meets $\pi$ in $B'$.
  Thus $A$ and $B'$ must lie in one of the 10 blocks
  of $D$ which is false.
  We conclude that either
  $$P_1A\cap P_2A'=P_3 \hbox{ and }P_1B\cap P_2B'=P_4$$
  or
  $$P_1A'\cap P_2A=P_3 \hbox{ and }P_1B'\cap P_2B=P_4.$$
  By examining the block $VCC'$ in a similar way
  we conclude that there are just
  2 possibilities for the 5-compressor, 
  namely $\mathcal{S}_1$ and $\mathcal{S}_2$
  as in
  Theorem~\ref{thm:ExistsTwo5Pts}, namely
  $$\mathcal{S}_1=\{P_1,P_2,P_1A\cap P_2A',P_1B\cap P_2B',P_1C\cap P_2C'\}$$
  or
  $$\mathcal{S}_2=\{P_1,P_2,P_1A'\cap P_2A,P_1B'\cap P_2B,P_1C'\cap P_2C\}.$$
\end{proofN}

%
%
%
\section{\ The number of planar Desargues configurations.\label{theNoPlanarDesConfigs}}
In this section $F=GF(q)$ the finite field of order $q$.
[The easiest example is when $q$ is a prime $p$ and our field
is formed from the integers modulo $p$ with the
usual addition and multiplication.]
First we count the number of Desargues configurations $D$
containing a given point $V$ of the plane $\pi$ contained
in $\Sigma=PG(3,q)$.
We have already shown the following.
\begin{itemize}
  \item[(a)]
    $D$ can be obtained as the section of a 5-compressor $\mathcal{S}$
    in $\Sigma$ with no point of $\mathcal{S}$ lying in $\pi$.
  \item[(b)]
    $\mathcal{S}=\{P_1,P_2,P_3,P_4,P_5\}$
    where $P_1,P_2$ are distinct fixed points 
    in space not in $\pi$
    such that the line $P_1P_2$ passes
    through the given point $V$ of $\mathcal{S}$.
  \item[(c)]
    There are exactly 2 such 5-compressors $\mathcal{S}_1,\mathcal{S}_2$
    whose section is the Desargues configuration $D$ in $\pi$ 
    containing $V$.
\end{itemize}

Our task now is to calculate the total number, say $\theta$,
of 5-compressors $\mathcal{S}$ containing $P_1,P_2$ as in (b) above.
It will then follow that there are exactly $\frac{\theta}{2}$
Desargues configurations in $\pi$ containing $V$.

In what follows, $\Lambda=AG(3,q)$ denotes
the affine 3-dimensional space over the field $F=GF(q)$
obtained from $\Sigma$ by removing $\pi$.

\begin{lemmaN}
  Let $P_1,P_2,P_3,P_4$ be a set of 4 given points in $\Lambda$
  not contained in a plane of $\Lambda$.
  Let $P_5$ be a point of $\Lambda$ such that
  $\mathcal{S}=\{P_1,P_2,P_3,P_4,P_5\}$
  is a 5-compressor in $\Lambda$, i.e., a set of 5 points in $\Lambda$
  such that no four of its points lie in a plane of $\Lambda$
  (and thus, such that no 3 points lie on a line of $\Lambda$).
  Then the numbers of possibilities for
  $P_5$ is 
  $(q-2)(q^2-2q+2)$.
  \label{noOfPossibilitiesForP5}
\end{lemmaN}
\begin{proofN}
  We calculate $u$, the number of points in the union $U$ of
  the 4 planes 
  ${<}P_1,P_2,P_3{>}$, ${<}P_1,P_2,P_4{>}$, ${<}P_1,P_3,P_4{>}$,
  ${<}P_2,P_3,P_4{>}$.
  Recall that \hfil\break ${<}P_i,P_j,P_k{>}$ denotes the unique plane
  of $\Lambda$ containing the points $P_i$, $P_j$ and $P_k$.

  No two of these 4 distinct planes are parallel
  so any two of them intersect in an affine line.
  Any 3 of the planes meet in a unique point in space which is one of
  the $P_i$.
  No point lies on all 4 planes.
  Then, by inclusion-exclusion we get that
  $u=|U|=4q^2-{4 \choose 2}q+{4 \choose 3}1=4q^2-6q+4$.
  Therefore the number of possibilities for $P_5$ is
  equal to $q^3-[4q^2-6q+4]=(q-2)(q^2-2q+2)$.
\end{proofN}

\begin{theoremN}
  The number of 5-compressors $\mathcal{S}=\{P_1,P_2,P_3,P_4,P_5\}$
  containing 2 given points $P_1,P_2$ in the affine space
  $\Lambda=AG(3,q)$
  is
  $$\theta=\frac{(q^3-q)(q^3-q^2)(q-2)(q^2-2q+2)}{6}.$$
  \label{no5PointsContainingP1P2}
\end{theoremN}
\begin{proofN}
  The point $P_3$ can be any point in $\Lambda$ not on the line
  $P_1P_2$.  So the number of choices for $P_3$
  is $q^3-q$.
  $P_4$ can be any point in $\Lambda$ not in the plane
  ${<}P_1,P_2,P_3{>}$.
  So there are $(q^3-q^2)$ choices for $P_4$.
  From Lemma~\ref{noOfPossibilitiesForP5},
  the number of choices for $P_5$ having chosen
  $P_1,P_2,P_3,P_4$ is $(q-2)(q^2-2q+2)$.
  The product $(q^3-q)(q^3-q^2)(q-2)(q^2-2q+2)$
  counts ordered triples $\{P_3,P_4,P_5\}$.
  Permuting these 3 points yields the same (unordered) set
  of three points.
  Since we are counting sets, not ordered sets,
  we divide this product by $3!=6$ to obtain the result.
\end{proofN}

\begin{theoremN}
  The total number of Desargues configurations
  in $PG(2,q)$ is
  $$\frac{q^3(q^3-1)(q^2-1)(q-2)(q^2-2q+2)}{120}.$$
  \label{theorem:noPlanarDesConfigs}
\end{theoremN}
\begin{proofN}
  Let $x$ denote the total number of Desargues configurations
  $D$ in $\pi=PG(2,q)$.
  We count $|J|$ the number of incidences in $J$ where
  $$J=\{(V,D)\ |\ \hbox{point }V\hbox{ is incident with Desargues configuration }D\}.$$
  For each point $V$ there are, from
  Theorem~\ref{no5PointsContainingP1P2},
  exactly $\frac{\theta}{2}$ configurations $D$ containing $V$.
  Thus 
  $|J|=(q^2+q+1)\frac{\theta}{2}$, since
  $\pi$ has $q^2+q+1$ points.
  Also $|J|=10x$ since $D$ contains exactly 10 points of $\pi$.
  Thus $x=\frac{(q^2+q+1)\theta}{20}$.
  Using the fact that $q^2+q+1=\frac{q^3-q}{q-1}$, and simplifying,
  the result follows.
\end{proofN}

%
%
%
\section{\ Three-dimensional Desargues configurations.\label{spatialDesConfigs}}
We will work here in $\Omega=PG(4,q)$, where
$\Sigma=PG(3,q)$ is a given 3-dimensional subspace of $\Omega$.
As in 
Section~\ref{sect:planarDesAnd5Pts},
a set $\mathcal{S}=\{P_1,P_2,P_3,P_4,P_5\}$
of 5 points in $\Omega$ with no 4 being coplanar
is called a 5-compressor.
The following lemma is analogous to 
Theorem~\ref{5Point}
in 3-dimensions.
\begin{lemmaN}
  Suppose no point of the 5-compressor $\mathcal{S}$
  lies in $\Sigma$.
  Then the 10 lines $P_iP_j$ joining $P_i$ to $P_j$ meet $\Sigma$
  in 10 distinct points $(ij)$.
  Each of the 10 planes ${<}P_i,P_j,P_k{>}$
  sections $\Sigma$ in a block
  $[ijk]$ of points $(ij)$, $(ik)$, $(jk)$
  which are collinear on a line called a
  {\em blockline}.
  The resulting configuration of 10 distinct points 
  and 10 distinct blocks
  is a Desargues configuration $D$ in $\Sigma$.
  \label{5PointIn3Dim}
\end{lemmaN}
As above we say that $\mathcal{S}$ sections $D$
or $\mathcal{S}$ sections to $D$.
$D$ is either ``planar'', i.e., lies in a plane in $\Sigma$,
or $D$ is ``non-planar'', so that 
the points of $D$ span a 3-dimensional space $\Sigma$.
The details are as follows.

\begin{lemmaN}
  \begin{itemize}
    \item[(a)]
      If $D$ is spatial, i.e., non-planar, 
      then $\mathcal{S}$ spans the space $\Omega=PG(4,q)$.
    \item[(b)]
      If $\mathcal{S}$ spans $\Omega$ then $D$ is spatial.
  \end{itemize}
  \label{spatialDesSSpans}
\end{lemmaN}
\begin{proofN}
  The lines $P_1P_2$, $P_1P_3$, $P_1P_4$, $P_1P_5$
  meet $\Sigma$ in the 4 points $(12)$, $(13)$, $(14)$
  and $(15)$.
  The points $(13)$, $(14)$, $(15)$
  are the vertices of one of the two triangles in $D$
  that are in perspective from $(12)$.
  Since $D$ is non-planar these 4 points,
  namely $(12)$, $(13)$, $(14)$, $(15)$ span $\Sigma$.
  $P_1$ is not in $\Sigma$.
  The line joining $P_1$ to $P_i$ contains $(1i)$, $i=2,3,4,5$.
  Thus $\mathcal{S}$ spans $\Omega={<}P_1,\Sigma{>}$.  
  This proves (a).

  To prove (b), 
  we suppose by way of contradiction that the points
  $(12)$, $(13)$, $(14)$ and $(15)$ lie in a plane $\pi$.
  Now $P_1$ is not in $\Sigma$.
  Consider the 3-dimensional space ${<}P_1,\pi{>}$.
  It contains $P_1$ and the points $(1i)$, $i=2,3,4,5$.
  Thus it contains $\mathcal{S}=\{P_1,P_2,P_3,P_4,P_5\}$.
  Then $\mathcal{S}$ is contained in a 3-dimensional space,
  a contradiction to the assumption that $\mathcal{S}$
  spans $\Omega$.
\end{proofN}

\begin{notationN}
  A 5-compressor $\mathcal{S}$ in $\Omega=PG(4,q)$
  that spans $\Omega$ is termed a
  {\em 5-arc} in $\Omega$.
  $\Gamma=AG(4,q)$ is
  the affine space obtained by removing $\Sigma$ from $\Omega$.
\end{notationN}
If $\mathcal{S}$ is a 5-arc then,
as for a 5-compressor, no 4 points of $\mathcal{S}$ are coplanar.
But, in addition, we demand that the 5 points of $\mathcal{S}$
are not contained in a 3-dimensional subspace of $\Omega$.

Using this notation 
and a similar method of proof to that of
Theorem~\ref{thm:ExistsTwo5Pts},
we have the following result.

\begin{lemmaN}
  Let $D$ be a spatial Desargues configuration in
  $\Sigma=PG(3,q)$ and let $V$ be a point of $D$.
  Let $P_1$, $P_2$ be fixed points in $\Omega=PG(4,q)$ 
  not in $\Sigma$
  such that the line $P_1P_2$ meets $\Sigma$ in $V$.
  Then
  \begin{itemize}
    \item[(a)] 
      there exists 5-compressors $\mathcal{S}_1$, $\mathcal{S}_2$
      in $\Gamma=AG(4,q)$
      such that $\mathcal{S}_1$, $\mathcal{S}_2$
      contain $P_1,P_2$
      and section $\Sigma$ in $D$.
    \item[(b)]
      Each of $\mathcal{S}_1$, $\mathcal{S}_2$ is a 5-arc in $\Gamma$.
    \item[(c)]
      Any 5-arc in $\Gamma$ containing $P_1,P_2$
      sections $\Sigma$ in a spatial Desargues configuration
      containing $V$.
  \end{itemize}
\end{lemmaN}

The following analogue of Theorem~\ref{thm:two5Points}
is proved in a similar way.
\begin{theoremN}
  Let $D$ be a spatial Desargues configuration in $\Sigma=PG(3,q)$
  where $\Sigma$ is contained in $\Omega=PG(4,q)$.
  Let $V$ be a point of $D$.
  Let $P_1, P_2$ be two points in the affine space
  $\Lambda=AG(4,q)$
  obtained by removing $\Sigma$ from $\Omega$ such that the
  line $P_1P_2$ meets $\Sigma$ in $V$.
  Then there are exactly two 5-compressors $\mathcal{S}_1,\mathcal{S}_2$
  in $\Omega$ containing $P_1$ and $P_2$ which section to $D$.
  Each of $\mathcal{S}_1,\mathcal{S}_2$ is a 5-arc
  of $\Omega$ which is contained in the affine space
  $\Gamma=AG(4,q)$.
\end{theoremN}

\begin{theoremN}
  The number of 5-arcs $\mathcal{S}$ of $\Omega$ contained in $\Gamma$
  and
  containing $P_1,P_2$ is
  $$\theta=\frac{(q^4-q)(q^4-q^2)(q^4-q^3)}{6}.$$
  \label{no5-arcsInAg4q}
\end{theoremN}
\begin{proofN}
  There are $q^4-q$ possibilities for $P_3$,
  followed by $q^4-q^2$ possibilities for $P_4$
  and $q^4-q^3$ possibilities for $P_5$.
  In total we have
  $(q^4-q)(q^4-q^2)(q^4-q^3)$ ordered possibilities
  for $\mathcal{S}$.
  Thus to find $\theta$ we divide by 
  $3!=6$ since this product counts ordered
  triples $P_3P_4P_5$.
\end{proofN}

\begin{theoremN}
  The total number of spatial Desargues configurations is
  $$\frac{(q^3+q^2+q+1)(q^4-q)(q^4-q^2)(q^4-q^3)}{120}.$$
  \label{noSpatialDes}
\end{theoremN}
\begin{proofN}
  Let $J$ denote the set of incidences
  $(V,D)$,
  where $V$ is a point of $\Sigma=PG(3,q)$ 
  and $D$ is a spatial Desargues configuration containing $V$.
  If there are $x$ spatial Desargues configurations
  in $PG(3,q)$ we get that
  $|J|=10x$.

  On the other hand, from
  Theorem~\ref{no5-arcsInAg4q}
  the number of spatial Desargues configurations containing $V$
  is $\frac{\theta}{2}$.
  Since there are $q^3+q^2+q+1$ possibilities for $V$ in $\Sigma$
  we get that $|J|=(q^3+q^2+q+1)\frac{\theta}{2}$.
  Thus
  $$x=\frac{(q^3+q^2+q+1)(q^4-q)(q^4-q^2)(q^4-q^3)}{120}.$$
\end{proofN}

We now sketch another proof of 
Theorem~\ref{noSpatialDes}.

\bigskip\bigskip
\noindent
{\em
  Alternative proof of 
  Theorem~\ref{noSpatialDes}.
}

\medskip\noindent
\begin{proofN}
  Our goal is to show
  that the number of spatial Desargues configurations
  containing $V$ is $\frac{\theta}{2}$.

  To this end let $VAA'$, $VBB'$ and $VCC'$
  be 3 blocks of a spatial Desargues configuration.
  There are $q^2+q+1$ choices for the line $VA$
  and $q^2+q$ choices for the line $VB$.
  Now the line $VC$ cannot lie in the plane formed
  by the lines $VA$ and $VB$ since the configuration
  is spatial.
  Thus there are $q^2$ choices for the line $VC$.
  All told we have
  $(q^2+q+1)(q^2+q)q^2$ choices for the
  (unordered) set of~3 lines.

  There are ${q \choose 2}$ choices for each of the sets
  $AA'$, $BB'$, $CC'$.  There are then 8 choices for the
  triangles $ABC$ and $A'B'C'$.  However,
  interchanging $A$ with $A'$, $B$ with $B'$,
  $C$ with $C'$ will yield the same configuration
  so we have effectively~4 choices for the triangles.
  Thus the number of spatial configurations
  containing $V$ is
  $$\frac{(q^2+q+1)(q^2+q)q^2{q\choose 2}{q\choose 2}{q\choose 2}\cdot 4}{6}$$
  which is equal to $\frac{\theta}{2}$, proving the result. 
\end{proofN}

%
%
%

\section{\ Desargues blockline structures \label{section:DesarguesBlockline}}

In this section, we address a historical ambiguity
as to the definition of a Desargues configuration
by both a Desargues configuration
and a Desargues blockline structure.

We are motivated by the following question:
Using the notation of Figure~\ref{fig:alternativeLabels},
given a Desargues configuration $D$ 
in $PG(2,F)$ with
\begin{itemize}
  \item[]
    points:\ \ \  $(12),(13),(14),(15),(23),(24),(25),(34),(35),(45)$
  \item[]
    and blocks: $\{(12)(13)(23)\}$, $\{(12)(14)(24)\}$, 
    $\{(12)(15)(25)\}$,\\
    \hphantom{Blocks of $D$: }
    $\{(13)(14)(34)\}$, $\{(23))(24)(34)\}$, \\
    \hphantom{Blocks of $D$: }
    $\{(13)(15)(35)\}$, $\{(23)(25)(35)\}$,\\
    \hphantom{Blocks of $D$: }
    $\{(14)(15)(45)\}$, $\{(24)(25)(45)\}$,\\
    \hphantom{Blocks of $D$: }
    $\{(34)(35)(45)\}$,
\end{itemize}
could there be a second Desargues configuration
that has the same points and the same blocklines as $D$,
but with different blocks than $D$?
A Desargues configuration's blocks each have exactly three points
by definition.  But, some blocklines could ``accidentally''
contain a fourth point of the Desargues configuration.  
This happens precisely when a point lies on its polar line.

Suppose a new block is formed by replacing one of the points of 
an existing block
with the fourth acciddental point.
The question is whether it is possible to find
a Desargues configuration
containing the new block which has the same points and 
blocklines as the original Desargues configuration?

For example, in Figure~\ref{fig:twoSCPtsBlockline},
the triangles $(13)(14)(15)$ and $(23)(24)(25)$
are in perspective from the point $(12)$.
Notice that the point $(35)$ is self-conjugate
as it is accidentally on the blockline that contains $[124]$.

There are four points of $D$ on the blockline
through $[124]$.  The given Desargues configuration
has a block $\{(12),(14),(24)\}$ from that line.
Could another Desargues configuration with the same points
and the same block lines have a block containing
$(35)$ and two of $(12),(14),(24)$?
We will show that the answer is no.

We introduce some terminology.
Associated with any Desargues configuration $D$ is
an incidence structure $D_{BL}$ of points and subsets of the points.
The points of $D_{BL}$ are the points of $D$.
The subsets of $D_{BL}$ are the points of $D$ that
lie on a blockline, i.e., a line of the plane containing a block.
(Such a subset has either three or four points.)
The resulting incidence structure is
called a {\em Desargues blockline structure}.
\hfil\break
Note: 
ostensibly, a better-sounding name for this might
be a Desargues blockline configuration.
However, as we have defined it in 
Section~\ref{section:DesarguesConfigurations},
a ``Desargues blockline configuration'' would not actually
be a configuration as 
the incidence structure of points and subsets of points
of a configuration must have
a constant number of points per subset.
In this situation, there will often be three points in a subset,
but there might be a fourth accidental point in a subset.

For example, if $D$ is the Desargues configuration
represented by Figure~\ref{fig:twoSCPtsBlockline},
the blocks of $D$ are 
    \hfil\break 
    \hphantom{Blocks of $D$: }
    $\{(12)(13)(23)\}$, $\{(12)(14)(24)\}$, 
    $\{(12)(15)(25)\}$,\\
    \hphantom{Blocks of $D$: }
    $\{(13)(14)(34)\}$, $\{(23))(24)(34)\}$, \\
    \hphantom{Blocks of $D$: }
    $\{(13)(15)(35)\}$, $\{(23)(25)(35)\}$,\\
    \hphantom{Blocks of $D$: }
    $\{(14)(15)(45)\}$, $\{(24)(25)(45)\}$,\\
    \hphantom{Blocks of $D$: }
    $\{(34)(35)(45)\}$.
\hfil\break
Each block has three points, as required by the definition
of a Desargues configuration.
If we focus on the block $[124]$ which contains
$(12),(14),(24)$, in $D$,
the corresponding subset in $D_{BL}$ contains four points,
namely $(12),(14),(24)$ and the accidental point $(35)$.

\begin{observationN}
  If a Desargues configuration has no self-conjugate points,
then the corresponding Desargues blockline structure
is the same as the given Desargues configuration.
\end{observationN}

The main result of this section is the following.
\begin{theoremN}
  Let $\pi=PG(2,F)$, where $F$ has characteristic
  different from $2,3$.
  Let $ABC$ and $A'B'C'$
  be two triangles in perspective from $V$
  that give rise to the Desargues configuration $D$.
  Assume that 
  two other triangles in perspective from $V$
  give rise to the Desargues configuration $D'$.
  Suppose that $D$ and $D'$ have the same points
  and the same blocklines.
  Then $D$ and $D'$ have the same blocks as well.
\label{theorem:DesarguesBlockline}
\end{theoremN}
\begin{proofN}
Consider a Desargues configuration $D$.
Since $F$ does not have characteristic~3,
$D$ can have at most three 
self-conjugate points~[\cite{bruenMcQGeomConfigs} Theorem 3.20].
We first show that there $D$ must have
at least one point that is not self-conjugate
and does not lie on any self-conjugate line.
(This is not necessarily the case
if the field has characteristic~3!)

To see this, there are four possibilities for the number of
self-conjugate points.  
We use coordinates as in Figure~\ref{fig:alternativeLabels}.
If there are no self-conjugate points, choose $V=(12)$
for the desired point.
If there is one self-conjugate point, label it $(13)$, say,
and then choose $V=(12)$ as the desired point.
If there are two self-conjugate points,
they must be on a block by~\cite[Lemma 3.10]{bruenMcQGeomConfigs}.
So, the labels for the two points must share a symbol.
Without loss of generality, they are $(13),(14)$.
In this case, choose $V=(12)$ as the desired point.
The remaining case is when there are three self-conjugate points.
Recall that each pair of them must be on a block.
So, each pair of them must share a symbol (in their labels).
But the three points are not on a block
since the characteristic is not~2 [\cite{bruenMcQGeomConfigs} Theorem 3.21].
Without loss of generality, they are $(13),(14),(15)$.
Again, we choose $V=(12)$ as the desired point.

Therefore, let $V=(12)$ be a point that 
is not self-conjugate and that
does not lie on any self-conjugate line.
There are exactly three blocklines of $D$ through $V$.
Now we consider $D'$, which has the same points and the
same blocklines as $D$.
As discussed in Section~\ref{section:DesarguesConfigurations},
$D'$ must have a pair of triangles in perspective
from $V$.  There are only three lines through $V$
and they only have two points in addition to $V$
on each of the three lines.
One blockline through $V$ contains the block $[123]$,
another contains $[124]$, and the final one contains $[125]$.

There are four possibilities for the possible pair
of triangles in $D'$:
\begin{itemize}
  \item[(i)]
    $(13)(14)(25),(23)(24)(15)$,
  \item[(ii)]
    $(13)(24)(15),(23)(14)(25)$,
  \item[(iii)]
    $(13)(24)(25),(23)(14)(15)$,
    or
  \item[(iv)]
    $(13)(14)(15),(23)(24)(25)$.
\end{itemize}

In case (i), there is a blockline in $D'$ through $(13)(25)$.
That means that $(13)$ or $(25)$ is self-conjugate (in $D$).
However, $(25)$ cannot be self-conjugate in $D$.
For this would imply that $(25)(13)(14)(34)$ lie on a
blockline in $D$ and thus on a blockline of $D'$.
But the points $(25)(14)(13)$ are postulated  to form
a triangle in $D'$.
Therefore $(13)$ is self-conjugate (in $D$).
Similarly, $(14)$ is self-conjugate (in $D$).
Also, there is a blockline through $(23)(15)$.
However, $(15)$ is not self-conjugate (in $D$) because
$(23)(24)(15)$ is a triangle.
Therefore, $(23)$ must be self-conjugate (in $D$).
But, $(14)$ and $(23)$ cannot both be self-conjugate (in $D$)
because they do not lie on a block (in $D$).
So, case (i) is not possible.

In case (ii),
there is a blockline through $(13)(24)$.
We know that $(24)$ is not self-conjugate (in $D$) because
$(13)(24)(15)$ is a triangle in $D'$.
Therefore $(13)$ is self-conjugate (in $D$).
Similarly, $(15)$ is self-conjugate (in $D$).
Also, there is a blockline through $(14)(23)$.
We know that $(14)$ is not self-conjugate (in $D$)
because $(23)(14)(25)$ is a triangle in $D'$.
Therefore $(23)$ is self-conjugate (in $D$).
But, $(23)$ and $(15)$ cannot both be self-conjugate (in $D$)
because they do not lie on a block.
So, case (ii) is not possible.

In case (iii),
there is a blockline through $(13)(24)$ and one through $(13)(23)$.
Therefore, $(23)$ and $(24)$ are self-conjugate (in $D$).
Also, there is a blockline through $(23)(14)$ and one
through $(23)(15)$.  Therefore $(14)$ and $(15)$
are self-conjugate (in $D$).  But, $(23)$ and $(14)$ cannot
both be self-conjugate (in $D$).
So, case (iii) is not possible.

That leaves case (iv), in which case $D=D'$.
\end{proofN}

\begin{corollaryN}
  Let $\pi=PG(2,F)$, where $F$ has characteristic
  different from $2,3$.
  If $D$ is a Desargues configuration that gives rise
  to the Desargues blockline structure $D_{BL}$,
  then no other Desargues configuration gives rise
  to the blockline structure $D_{BL}$.
  \label{corollary:uniqueDesConfigfromBL}
\end{corollaryN}

It can be shown by similar synthetic methods
that the restriction in Theorem~\ref{theorem:DesarguesBlockline}
on the characteristic of the field
is not really necessary.
Corollary~\ref{corollary:uniqueDesConfigfromBL}
then provides
a significant strenghtening to
our calculation of the total number of Desargues configurations
containing a point in Section~\ref{theNoPlanarDesConfigs}.

\begin{corollaryN}
  Let $\pi=PG(2,F)$, where $F$ is a finite field.
  Then the total number of distinct
  Desargues blockline structures 
  is the same as the total number
  of Desargues configurations.
  \label{corollary:noDesBlockline}
\end{corollaryN}
  Added in proof.
  Using similar methods we can show that
  Corollary~\ref{corollary:noDesBlockline} holds in all cases
  regardless of the characteristic.

%
%
%

\section{\ A possibly easier approach to the enumeration of the number of planar Desargues configurations through a given point that was given in Section~\ref{theNoPlanarDesConfigs}.\label{section:easierApproach}}
In this section, we consider the possibility of an easier calculation
of the number of
planar Desargues configurations through a given point
that was given in Section~\ref{theNoPlanarDesConfigs}.

We have seen that the answer is $\theta/2$,
where
  $$\theta=\frac{(q^3-q)(q^3-q^2)(q-2)(q^2-2q+2)}{6}$$
as in Theorem~\ref{no5PointsContainingP1P2}
and Theorem~\ref{theorem:noPlanarDesConfigs}.
Is there a quicker way to calculate this
by just 
working in the plane?

The following calculation seems plausible
at first glance, but it is incorrect.

ALTERNATE {\bf (INCORRECT)} CALCULATION:

Fix a point $V$ in the plane $PG(2,q)$.
Consider a prospective Desargues configuration $D$ containing that point.
We wish to choose the points of the triangles $ABC$, $A'B'C'$
that are in perspective from $V$, with $V=AA'\cap BB'\cap CC'$.

We first choose the three lines throught $V$
that will be $AA'$, $BB'$ and $CC'$.
There are ${q+1 \choose 3}$ choices for those lines.

On one of those three lines we choose two points
unequal to $V$
to be the set $\{A, A'\}$.
This gives us ${q \choose 2}$ choices.
On the second line we choose the set
$\{B,B'\}$ giving ${q\choose 2}$ choices.
Then any point of $\{A,A'\}$ can be paired with
any one of $\{B,B'\}$ - there are two choices here - to yield the sides
$AB,A'B'$.

Finally, we need to choose $C$ and $C'$ on the
remaining line.
Let us choose $C$ first.
We ensure that $C$ is not $V$
and $C$ is not on the line $AB$.
That gives us $q+1-2=q-1$ choices.

To choose $C'$, we need to make sure that
$C'\ne V, C$ and that $C'$ is not on the line $A'B'$.
That gives us $q+1-3=q-2$ choices.

This yields a total of
$${q+1 \choose 3}{q \choose 2}{q \choose 2}2(q-1)(q-2).$$
However, this is not the answer obtained in
Section~\ref{theNoPlanarDesConfigs}.

What went wrong?  {\em We did not pay attention to the possibility
that a point might be a self-conjugate point!}
For example, when we chose $C'$, we did not consider the possibility
that $C$ might be self-conjugate.  If it is self-conjugate,
then $C$ would equal $CC'\cap A'B'$.  So, it was over-simplistic
to just subtract three in our calculation for the number
of choices for $C'$.
Furthermore, what if $B$ happened to be a self-conjugate point?
Note that all of $A,B,C$ might perhaps be self-conjugate or not.
To make things even more interesting, characteristic~3 fields
play a special role.
In other words the characteristic of the field
has to be taken into account using
the na\"ive approach above!
It is possible to have four self-conjugate points
if and only if the field has characteristic~3.
This was shown in~\cite[Theorem 3.20]{bruenMcQGeomConfigs}.

%
%
%

\section{\ Concluding remarks.\label{section:concludingRemarks}}
Desargues theorem in the plane can easily be shown
using coordinates as in Lord (\cite{lord}).
The synthetic proofs seem to all use a construction which
is in effect a compressor or a 5-point.
Baker (\cite{baker}) uses the Desargues theorem in three dimensions
to finish the proof.
The Veblen-Young proof (\cite{veblenYoungVol1})
is a very slightly longer version of our proof here.
But the important continuation here,
using the pairs-triples notation is the
connection with polarities developed 
in~\cite{bruenMcQGeomConfigs,bruenMcQFourSC}
and used in the enumeration.
One possible method is to count the configurations having
0 self-polar points or exactly 1, 2, 3 or 4 self-polar points.
From the Puystjens-Thas (\cite{thas}) paper we 
have seen how difficult this is.
The fact that 4 self-conjugate points only occur in
characteristic~3 shows that for general fields
of order~$q$ the method is doomed to failure.
Fortunately the necessary information for enumeration
is stored in compressed form
in the 5-compressor and we can count those.

In~\cite{thas} the authors concentrate on 
the cases when the underlying field has characteristic~2 or~3.
They also demand that, in our terminology,
each blockline must contain exactly 3 points of
the Desargues configuration.

In the general case it is possible that a blockline
contains as many as 4 points of the configuration, consisting
of a block and an ``accidental'' extra point.

\begin{figure}[!ht]
    \center
    \includegraphics[height=9cm]{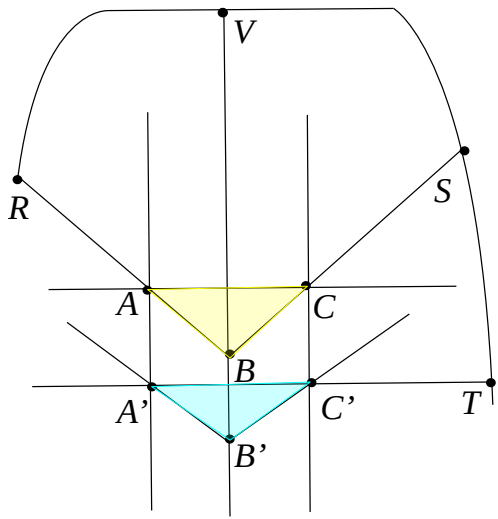}
    \caption{
    The Desargues configuration containing the triangles
    $ABC$, $A'B'C'$ has a blockline with four points,
    namely $R,S,T,V$.
    The triangle $A'B'C'$ is a translation of $ABC$.
    Their Desargues axis contains $V$.
    \label{fig:transl}}
\end{figure}

An easy example in which there
are four points of a Desargues configuration
on a blockline
is 
as follows.
Let $T$ be a translation with centre $V$ and axis $l_{\infty}$,
where $l_{\infty}$ is the line at inifinity
and $V$ is on $l_{\infty}$.
Let $ABC$ denote an affine triangle
and let $A'B'C'$ be its image under $T$.
These two triangles are in perspective from $V$.
The Desargues axis is $l_{\infty}$ and it
contains 4 points of the resulting configuration consisting of the~3 intersections
of corresponding sides together with the
``accidental extra point'' $V$.
See Figure~\ref{fig:transl}.
What is happening here is that the point
$V$ is self-polar or self-conjugate
(see Section~\ref{section:DesPlanePolarities}).
A detailed discussion of self-conjugate points 
in a Desargues configuration
is provided in~\cite{bruenMcQGeomConfigs}.

The case of characteristic~3 actually plays a 
very special role.
In~\cite{bruenMcQGeomConfigs} we show that
there can be at most three self-conjugate points
in a Desargues configuration
unless the characteristic is 3,
in which case there can exist as many
as four self-conjugate points.

The number of configurations in~\cite{thas}
is smaller than the number in
Theorem~\ref{theorem:noPlanarDesConfigs}
stemming from the fact that the authors
there are working with a restricted class
of configurations.

%
%

\end{document}